
\magnification=1200
\hsize=12.5cm
\vsize=18cm
\hoffset=1cm
\voffset=2cm

\nopagenumbers               
\footnote{{}}{\hskip -2 mm 2000 {\it Mathematics Subject
Classification}. Primary 11M06, Secondary 11F72.}

\def\DJ{\leavevmode\setbox0=\hbox{D}\kern0pt\rlap
 {\kern.04em\raise.188\ht0\hbox{-}}D}

\def\txt#1{{\textstyle{#1}}}
\baselineskip=13pt
\def\hf{{\textstyle{1\over2}}}
\def\a{\alpha}
\def\d{{\,\rm d}}
\def\e{\varepsilon}
\def\f{\varphi}

\def\k{\kappa}
\def\s{\sigma}

\def\={\;=\;}
\def\zx{\zeta(\hf+ix)}
\def\zt{\zeta(\hf+it)}

\def\R{\Re{\rm e}\,} \def\I{\Im{\rm m}\,} \def\s{\sigma}
\def\z{\zeta}

\def\no{\noindent} 
\def\H{H_j^3({\txt{1\over2}})}  \def\={\,=\,}
\def\hf{{\textstyle{1\over2}}}
\def\txt#1{{\textstyle{#1}}}
\def\f{\varphi}
\def\Z{{\cal Z}}
\def\K{{\cal K}}
\font\tenmsb=msbm10
\font\sevenmsb=msbm7
\font\fivemsb=msbm5
\newfam\msbfam
\textfont\msbfam=\tenmsb
\scriptfont\msbfam=\sevenmsb
\scriptscriptfont\msbfam=\fivemsb
\def\Bbb#1{{\fam\msbfam #1}}

\def \NN {\Bbb N}
\def \CC {\Bbb C}
\def \RR {\Bbb R}
\def \ZZ {\Bbb Z}

\font\aa=cmcsc10 at 12pt 
\font\bb=cmcsc10
\font\cc=cmcsc10 at 8pt 

\font\ff=cmr9
\font\hh=cmbx12 
\def\rightheadline{\hfil\ff The estimation of Mellin transforms\hfil\folio}
\def\leftheadline{\ff\folio\hfil Aleksandar Ivi\'c
 \hfil}
\def\emptyheadline{\hfil}
\headline{\ifnum\pageno=1 \emptyheadline\else
\ifodd\pageno \rightheadline \else \leftheadline\fi\fi}
\topglue2cm
\centerline{\aa  On the estimation of some Mellin transforms}

\medskip
\centerline{\aa
connected with the fourth moment of $|\zt|$}
\bigskip\bigskip
\centerline{\aa Aleksandar Ivi\'c}
\bigskip
\centerline{\bb Abstract}
\bigskip\bigskip\no
{\ff Mean square estimates for $\Z_2(s) = \int_1^\infty|\zx|^4x^{-s}\d x\;(\R s > 1)$  
are discussed, and some related Mellin transforms of quantities
connected with the fourth power moment of $|\zx|$.}

\bigskip\bigskip
{\hh 1. Introduction}
\bigskip

\noindent

Let $\Z_2(s)$  be the analytic continuation of the function
$$
\Z_2(s) \= \int_1^\infty |\zx|^4x^{-s}\d x \qquad(\R s > 1),
$$
which represents the (modified) Mellin transform of $|\zx|^4$.
It was introduced by Y. Motohashi
[14] (see also  [7], [9], [12]  and [15]), who showed  that  it  
has meromorphic continuation over $\CC$. 
In the half-plane $\s = \Re{\rm e}\, s > 0$
it has the following singularities: the pole $s = 1$ of order five, 
simple poles at $s = {1\over2} \pm i\k_j\,(\k_j
= \sqrt{\lambda_j - {1\over4}})$ and poles at $s = \hf\rho$, 
where $\rho$ denotes complex
zeros of $\zeta(s)$.  Here as usual $\,\{\lambda_j = \k_j^2 +
{1\over4}\} \,\cup\, \{0\}\,$ is the discrete spectrum of the
non-Euclidean Laplacian acting on $SL(2,\ZZ)$-automorphic forms
(see [15, Chapters 1--3] for a 
comprehensive account of  spectral theory and the Hecke $L$-functions).

\medskip
The  aim of this note is to study the estimation $\Z_2(s)$
 in mean square and the (modified) Mellin transforms of certain
other quantities related to the fourth power moment
of $|\zt|$. This  research was begun in [12], and 
continued in [7] and [9]. It was shown there that we have
$$
\int_0^T|{\cal Z}_2(\s + it)|^2\d t \ll_\e 
T^\e\left(T + T^{2-2\s\over1-c}\right) \qquad(\hf < \s < 1),\leqno(1.1)
$$
and we also have unconditionally
$$
\int_0^T|\Z_2(\s + it)|^2\d t \;\ll\;T^{10-8\s\over3}\log^CT
\qquad(\hf < \s <1,\,C > 0).\leqno(1.2)
$$
Here and later $\e$ denotes arbitrarily small, positive constants,
which are not necessarily the same ones at each occurrence, while
$\s$ is assumed to be fixed. The constant $c$ 
appearing in (1.1) is defined by
$$
E_2(T) \ll_\e T^{c+\e},\leqno(1.3)
$$
where the function $E_2(T)$ denotes the error term in the asymptotic
formula for the mean fourth power of $|\zt|$. It is customarily
defined by the relation
$$
\int_0^T|\zt|^4\d t \;=\; TP_4(\log T) \;+\;E_2(T),\leqno(1.4)
$$
with
$$
P_4(x) \;=\; \sum_{j=0}^4\,a_jx^j, \quad a_4 = {1\over2\pi^2}.\leqno(1.5)
$$
For an explicit evaluation of the $a_j$'s in (1.5), 
see the author's work [3]. The best known value of $c$ in (1.3) is
$c = 2/3$ (see e.g., [11] or [15]), and it is conjectured that
$c = 1/2$ holds, which would be optimal. Namely 
(see [2], [4], [14] and [15]) one has
$$
E_2(T) = \Omega_\pm(T^{1/2}).\leqno(1.6)
$$

Mean value estimates for $\Z_2(s)$ are a natural tool to investigate
the eighth power moment of $|\zt|$. Indeed, one has (see [7, eq. (4.7)])
$$
\int_T^{2T}|\zt|^8\d t \ll_\e 
T^{2\s-1}\int_1^{T^{1+\e}}|\Z_2(\s+it)|^2\d t \quad(\hf < \s < 1).
\leqno(1.7)
$$ 
In [9] the pointwise estimate for $\Z_2(s)$ was given  by
$$
\Z_2(\s + it) \;\ll_\e\; t^{{4\over3}(1-\s)+\e},\leqno(1.8)
$$
for $\hf < \s \le 1$ fixed and $t \ge t_0 > 0$.
This result is still much weaker than the bound conjectured in [7] by the
author, namely that
for any given $\e > 0$ and fixed $\s$ satisfying $\hf < \s < 1$, one has
$$
\Z_2(\s + it) \;\ll_\e\; t^{{1\over2}-\s+\e}\qquad(t \ge t_0 > 0).
$$
\bigskip
To define another Mellin transform related to $\Z_2(s)$,
let $P_4(x)$ be defined by (1.4),  let 
$$
Q_4(x) := P_4(x) + P_4'(x)\leqno(1.9)
$$ 
and set
$$
\K(s) := \int_1^\infty (|\zx|^4 - Q_4(\log x))E_2(x)x^{-s}\d x. \leqno(1.10)
$$
The integral on the right-hand side of (1.10) converges absolutely at
least for $\s > 5/3$, in view of (1.3) with $c=2/3$ and the bound for the
fourth moment of $|\zx|$. However, the interest in $\K(s)$ lies in the
fact that (1.4) and (1.9) yield
$$
E_2'(x) = |\zx|^4 - Q_4(\log x).\leqno(1.11)
$$
Thus an integration by parts shows that
$$
\K(s) = -\hf E_2^2(1) - \hf s\int_1^\infty E_2^2(x)x^{-s-1}\d x.\leqno(1.12)
$$
In view of the mean square bound (see e.g., [10] and [15])
$$
\int_1^T E_2^2(t)\d t \ll T^2\log^CT,\leqno(1.13)
$$
it follows that (1.12) furnishes analytic continuation of $\K(s)$ to the
region $\s > 1$. A true asymptotic formula for the integtal in (1.13)
would provide further analytic continuation of $\K(s)$. For example,
a strong conjecture is that
$$
\int_1^T E_2^2(t)\d t = T^2p(\log T) + R(T),\quad R(T) \ll_\e T^{\rho+\e}
\leqno(1.14)
$$
with $p(x)$ a suitable polynomial (perhaps of degree zero) and
${3\over2} \le \rho < 2$. Namely from [2, Theorem 4.1] with $|H| =
T^{\rho/3}$ it follows that (1.14) implies (1.3) with $c \le \rho/3$, hence
$\rho \ge {3\over2}$ must hold in view of (1.7). Then the integral in (1.12)
becomes
$$
\int_1^\infty (2p(\log x) + p'(\log x))x^{-s}\d x
+ O(1) + (s+1)\int_1^\infty R(x)x^{-s-1}\d x.
$$
The first integral above is easily evaluated as 
$\sum_{j=1}^{m+1}b_j(s-1)^{-j}$, where $m$ is the degree of $p(x)$. 
The second integral is regular for  $\s > \rho-1$ if 
$R(x) \ll_\e x^{\rho+\e}$. Thus, on (1.14), it is seen 
that $\K(s)$ is regular
for $\s > \rho-1$ except for a pole at $s=1$ of order $1 + {\rm{deg}}\,p(x)$.

\medskip
Finally we define a Mellin transform related to the spectral theory
of the non-Euclidean Laplacian. Let, as usual,
$\a_j = |\rho_j(1)|^2(\cosh\pi\k_j)^{-1}$, where
$\rho_j(1)$ is the first Fourier coefficient of the Maass wave form
corresponding to the eigenvalue $\lambda_j$ to which the Hecke series
$H_j(s) = \sum_{n=1}^\infty t_j(n)n^{-s}$ is attached. For $0 < \xi \le 1$
we define
$$
{\cal I}(t;\xi) = {1\over\sqrt{\pi}t^\xi}\int_{-\infty}^\infty
|\z(\hf + it + iu)|^4\exp(-(u/t^\xi)^2)\d u.\leqno(1.15)
$$
The importance of the function ${\cal I}(t;\xi)$ in the theory of the fourth
power moment of $|\zt|$ comes from the fact that, for $\hf \le \xi < 1$
and suitable $C>0$, we have (see Y. Motohashi [13] and [15])
the explicit formula
$$\eqalign{
{\cal I}(t;\xi) &=
 {\pi\over\sqrt{2t}}\sum_{j=1}^\infty \a_j H_j^3(\hf)\k_j^{-1/2}
\sin\left(\k_j\log{\k_j\over4{\rm e}t}\right)
\exp(-{\txt{1\over4}}(t^{\xi-1}\k_j)^2) + O(\log^Ct)\cr&
= I(t;\xi) + O(\log^Ct),\cr}\leqno(1.16)
$$
say. Then we let
$$
J(s,\xi) = \int_1^\infty I(x;\xi)x^{-s}\d x\leqno(1.17)
$$
denote the (modified) Mellin transform of $I(t;\xi)$.
In view of the bound (see e.g. [15])
$$
\sum_{\k_j\le T}\a_j\H \ll T^2\log^CT\qquad(C>0)\leqno(1.18)
$$
it easily follows that $J(s,\xi)$ is regular for $\s > 2 - {3\over2}\xi$.

\medskip
The plan of the paper is as follows. The function $\Z_2(s)$ will be studied
in Section 2, $\K(s)$ in Section 3, while Section 4 is devoted to
$J(s,\xi)$.

\bigskip
{\hh 2. The function $\Z_2(s)$}
\bigskip
\no
The new result concerning mean square bounds for $\Z_2(s)$ is contained in

\bigskip
THEOREM 1. {\it For ${5\over6} \le \s \le {5\over4}$ we have}
$$
\int_1^{T}|\Z_2(\s+it)|^2\d t \ll_\e T^{{15-12\s\over5}+\e}.\leqno(2.1)
$$

\medskip 
{\bf Proof}. To prove (2.1) we first introduce, as in  [7], the function
$$
F_K(s) \;:=\; \int_{K/2}^{5K'/2}\f(x)|\zx|^4x^{-s}\d x
\qquad(K < K' \le 2K),\leqno(2.2)
$$
where $\f(x) \in C^\infty$ is  a nonnegative function supported
in $[K/2,\,5K'/2]$ such that $\f(x) = 1$ for $K < K' \le 2K$, and
$$
\f^{(r)}(x) \;\ll_r\;K^{-r}\qquad(r = 0,1,2,\ldots).\leqno(2.3)
$$
To connect $F_K(s)$ and $\Z_2(s)$ note that from
the Mellin  inversion formula (e.g., [7, eq. (2.6)] we have
$$
|\zx|^4 \= {1\over2\pi i}\int_{(1+\e)}{\cal Z}_2(s)x^{s-1}\d s \qquad(x>1),
$$
where the $\int_{(c)}$ denotes integration over the line $\R s = c$.
Here we replace the line of integration  by the contour ${\cal L}$,
consisting of the same straight line from which the segment
$[1+\e-i,\,1+\e+i]$ is removed and replaced by a circular arc 
of unit radius, lying
to the left of the line, which passes over the pole $s =1 $ of
the integrand. By the residue theorem we have
$$
|\zx|^4 \= {1\over2\pi i}\int_{\cal L}{\cal Z}_2(s)x^{s-1}\d s
+ Q_4(\log x) \qquad(x > 1),\leqno(2.4)
$$
where   $Q_4$ is defined by (1.9). Hence by using (2.4) we obtain
$$\eqalign{
F_K(s) &= {1\over2\pi i}\int_{\cal L}{\cal Z}_2(w)
\left(\int_{K/2}^{5K'/2}\f(x)x^{w-s-1}\d x\right)\d w\cr&
+ \int_{K/2}^{5K'/2}\f(x)Q_4(\log x)x^{-s}\d x.\cr}\leqno(2.5)
$$
In view of (2.3) we infer, by repeated integration by parts,
that the last integral in (2.5) is $\ll T^{-A}$ for any given $A>0$.
Similarly we note that
$$\eqalign{
&\int_{K/2}^{5K'/2}\f(x)x^{w-s-1}\d x\cr&
= (-1)^r\int_{K/2}^{5K'/2}\f^{(r)}(x){x^{w-s+r-1}\over(w-s)\cdots
(w-s+r-1)}\d x \ll T^{-A}\cr}
$$
for any given $A>0$, provided that $|\I w - \I s| > T^\e$ and
$r = r(A,\e)$ is sufficiently large. Thus if in the $w$--integral
in (2.5) we replace the contour ${\cal L}$ by the straight line
$\R w = d\;(\hf < d < 1)$, we shall obtain
$$
F_K(s) \;\ll_\e\; K^{d-\s}\int_{t-T^\e}^{t+T^\e}|\Z_2(d+iv)|\d v + T^{-2}.
\leqno(2.6)
$$
Squaring (2.6) and integrating it follows that
$$
\int_T^{2T}|F_K(s)|^2\d t \ll_\e T^{-1} + K^{2d-2\s}T^\e
\int_{T-T^\e}^{2T+T^\e}|\Z_2(d+iv)|^2\d v\quad(\s>d).\leqno(2.7)
$$
However, from (1.1) and (1.3) with $c = 2/3$  it follows that
$$
\int_T^{2T}|\Z_2({\txt{5\over6}} + it)|^2\d t \;\ll_\e\;
T^{1+\e},\leqno(2.8)
$$
so that (2.7) yields
$$
\int_T^{2T}|F_K(s)|^2\d t \ll_\e T^{-1} + K^{5/3-2\s}T^{1+\e}
\qquad({\txt{5\over6}} \le \s \le 1).\leqno(2.9)
$$
Now we write
$$
\Z_2(s) = \int_1^X |\zx|^4x^{-s}\d x 
+ \int_X^\infty |\zx|^4x^{-s}\d x = I_1 + I_2,
$$
say, which initially holds for $\s>1$. To estimate the mean square
of $I_1$, we use the bound (which, up to `$\e$', is the strongest one
known; see e.g., [1])
$$
\int_1^X |\zx|^8\d x  \ll_\e X^{3/2+\e},\leqno(2.10)
$$
and the following lemma, whose proof can be found in [7].

\medskip
LEMMA 1. {\it Suppose that $g(x)$ is a real-valued, 
integrable function on $[a,b]$, a subinterval
of $[2,\,\infty)$, which is not necessarily finite. Then}
$$
\int\limits_0^{T}\left|\int\limits_a^b g(x)x^{-s}\d x\right|^2\d t
\le 2\pi\int\limits_a^b g^2(x)x^{1-2\s}\d x \quad(s = \s+it\,,T > 0,\,a<b).
\leqno(2.11)
$$
\medskip
Then, from (2.10) and (2.11), we find that
$$
\int_T^{2T}I_1^2\d t \ll_\e X^{5/2-2\s+\e}.\leqno(2.12)
$$
From (2.6) and (2.7) with $d = 5/6$ we obtain the analytic continuation of
$I_2 = I_2(s)$ to the region $\s > 5/6$, taking first $K = 2X$, 
writing $1 = \f(x) + (1-\f(x)$ in the integral over
$[\hf X,\,X]\,$, and estimating the mean square of
$$
\int_{X/2}^X(1-\f(x))|\zx|^4x^{-s}\d x
$$
by the bound in (2.11). For the remaining integrals we use, after
the integrals are split into subintegrals of the type $F_K(s)$, 
the bound given by (2.9). We obtain
$$
\int_T^{2T}|\Z_2(\s+it)|^2\d t \ll_\e X^{5/2-2\s+\e}
+ T^{1+\e}X^{5/3-2\s} \ll_\e T^{(15-12\s)/5+\e}\leqno(2.13)
$$ 
with the choice $X = T^{6/5}$. Replacing $T$ by $T2^{-j}$ and adding
up all the results, we obtain (2.1) in the range ${5\over6} \le \s \le 1$.

\medskip
To obtain (2.1) in the remaining range $1 < \s \le {5\over4}$, first
we note that by a slight change of proof we see that (2.7) holds
for $d\ge 1$. Thus invoking (2.1) with $\s = 1$ it is seen that
for $1 < \s \le {5\over4}$ (when the exponent in (2.12) is non-negative)
we obtain
$$
\int_T^{2T}|\Z_2(\s+it)|^2\d t \ll_\e X^{5/2-2\s+\e}
+ T^{3/5+\e}X^{2-2\s} \ll_\e T^{(15-12\s)/5+\e}\leqno(2.14)
$$ 
again with the choice $X = T^{6/5}$. The proof of Theorem 1 is complete.

\medskip
As a corollary of (2.1) we can obtain (2.10), although this is somewhat
going round in a circle, since we actually used (2.10) in the course of
proof of (2.1). Recall that we have (1.7), 
but the analysis of its proof clearly shows that it remains valid for
$\s \ge 1$ as well. If we use (2.1) with $\s = 5/4$ in (1.7), then
(2.10) immediately follows. The essentialy new result provided by Theorem 1
is the bound
$$
\int_1^T|\Z_2(1+it)|^2\d t \ll_\e T^{3/5+\e},\leqno(2.15)
$$
and it would be of great interest to decrease the exponent of $T$ on
the right-hand side of (2.15). In fact, the hypothetical estimate
$$
\int_1^X |\zx|^8\d x  \ll_\e X^{1+\e}\leqno(2.16)
$$
is equivalent to
$$
\int_1^T|\Z_2(1+it)|^2\d t \ll_\e T^\e.\leqno(2.17)
$$
From (1.7) with $\s=1$ it follows at once that (2.17) implies (2.16),
and the other implication follows by the method of proof of (1.7) 
in [7]. This fact stresses out once again the importance of 
mean square bounds for $\Z_2(s)$.

\bigskip
{\hh 3. The function $\K(s)$}
\bigskip
In this section we shall deal with the function $\K(s)$, defined
by (1.10) or (1.12). Our result is the following

\bigskip
THEOREM 2. {\it The function $\K(s)$, defined by} (1.10), {\it admits
analytic continuation which is regular for $\R s >1$. It satisfies}
$$
\K(\s+it) \ll_\e |t|^\e(|t|^{3-2\s}+1)\qquad(\s>1) \leqno(3.1)
$$
{\it and}
$$
\int_0^T|\K(\s + it)|^2\d t \ll_\e T^{{13-6\s\over3}+\e}
\qquad({\txt{7\over6}} \le\s\le {\txt{13\over6}}). \leqno(3.2)
$$

\bigskip {\bf Proof}.
To prove (3.1) note first that, by the Cauchy-Schwarz inequality for 
integrals, we have ($C>0$)
$$
\eqalign{&
\int_Y^{2Y}\left||\zx|^4-Q_4(\log x)\right||E_2(x)|x^{-\s}\d x
\ll Y^{3/2-\s}\log^CY +
\cr& \,+ Y^{-\s}\left(\int_Y^{2Y}|\zx|^4\d x\right)^{1/2}
\left(\int_Y^{2Y}(|\zx|^4-Q_4(\log x))E_2^2(x)\d x\right)^{1/2}\cr&
\ll Y^{3/2-\s}\log^CY.\cr}\leqno(3.3)
$$
In the last integral we integrated by parts,  recalling that that (1.11) 
holds, as well as (1.13) and (1.3) with $c=2/3$. The above bound shows
then that $\K(s) \ll 1$ for $\s > 3/2$. Suppose now that 
$1 < \s \le 3/2$. Similarly to (1.12) we have
$$\eqalign{
\K(s) &= \int_1^X(|\zx|^4-Q_4(\log x))E_2(x)x^{-s}\d x\cr&
- \hf E_2^2(X)X^{-s} + \hf s \int_X^\infty E_2^2(x)x^{-s-1}\d x.\cr}
\leqno(3.4)
$$
From (3.3) it follows that the first integral above is $\ll 
X^{3/2-\s}\log^CX$, and the second (by (1.13)) is $\ll |t|X^{1-\s}$.
The choice $X=t^2$ easily leads then to (3.1).

\medskip
To prove (3.2) we start from (3.4) and use Lemma 1. We obtain
$$\eqalign{&
\int_T^{2T}\left|
\int_1^X(|\zx|^4-Q_4(\log x))E_2(x)x^{-s}\d x\right|^2\d t\cr&
\ll \int_1^X(|\zx|^4-Q_4(\log x))^2E^2_2(x)x^{1-2\s}\d x.\cr}\leqno(3.5)
$$
Defining the Lindel\"of function
$$
\mu(\s) \,=\,\limsup_{t\to\infty}\,{\log|\z(\s+it)|\over\log t}
\qquad(\s\in\RR)
$$
in the customary way and letting $\f(x)$ be as in (2.2), we see that
$$\eqalign{&
\int_{K/2}^{5K'/2}\f(x)(|\zx|^4-Q_4(\log x))^2E^2_2(x)x^{1-2\s}\d x\cr&
\ll_\e K^{1-2\s+4\mu({1\over2})+\e}
\int_{K/2}^{5K'/2}\f(x)(|\zx|^4 + \log^8x)E^2_2(x)\d x\cr&
= K^{1-2\s+4\mu({1\over2})+\e}\int_{K/2}^{5K'/2}\f(x)
(E_2'(x) + Q_4(\log x) + \log^8x)E^2_2(x)\d x\cr&
= K^{1-2\s+4\mu({1\over2})+\e}\int_{K/2}^{5K'/2}
\f(x)({\txt{1\over3}}E^3_2(x))'
\d x + O_\e(K^{3-2\s+4\mu({1\over2})+\e})\cr&
\ll_\e K^{3-2\s+4\mu({1\over2})+\e}.
\cr}
$$
Therefore the expression on the right-hand side of (3.5) is, for
$X = T^C, C>0$,
$$
\ll_\e T^\e(1 + X^{3-2\s+4\mu({1\over2})}).\leqno(3.6)
$$
Next we have, by Lemma 1, (1.3) with $c = 2/3$ and (1.13),
$$\eqalign{&
\int_{T}^{2T}\left|s\int_X^\infty E_2^2(x)x^{-s-1}\d x\right|^2\d t
\ll T^2\int_X^\infty E_2^4(x)x^{-2\s-1}\d x
\cr&
\ll_\e T^2\int_X^\infty E_2^2(x)x^{1/3-2\s+\e}\d x
\ll_\e T^2X^{7/3-2\s+\e},
\cr}\leqno(3.7)
$$
provided that $\s > {7\over6}$. From (3.4), (3.6) and (3.7) we infer that
$$\eqalign{&
\int_{T}^{2T}|\K(\s+it)|^2\d t
\ll_\e T^\e(1+ X^{3-2\s+4\mu({1\over2})}+T^2X^{7/3-2\s})\cr&
\ll_\e T^{{13-6\s\over3}+\e}\qquad({\txt{7\over6}} \le \s \le
{\txt{13\over6}})
\cr}
$$
with the trivial bound $\mu(\hf) < {1\over6}$ and $X = T$. This easily gives
(3.2), and slight improvements are possible with a better value of
$\mu(\hf)$. A mean square bound can also be obtained for the whole
range $\s > 1$, by using the trivial bound $tX^{1-\s+\e}$ for
the second integral in (3.4). This will lead to
$$
\int_1^T|{\cal K}(\s+it)|^2\d t \ll_\e \cases{T^{1+\e}&$(\s > 3/2),$
\cr{}
\cr
T^{{33-18\s\over5}+\e}&$(1 < \s \le 3/2)$.\cr}
$$

\bigskip
Mean square estimates for $\K(s)$ can be used to bound the fourth
moment of $E_2(t)$, much in the same way that mean square estimates
for $\Z_2(s)$ can be used (cf. (1.7)) to bound the eighth moment
of $|\zt|$. We have

\bigskip
THEOREM 3. {\it For $\s>1$ fixed}
$$
\int_T^{2T} E_2^4(t)\d t  \;\ll_\e\; T^{2\s+1}
\left(1 + \int_0^{T^{1+\e}}\,{|\K(\s+it)|^2\over1+t^2}\d t\right).
\leqno(3.8)
$$

\bigskip
{\bf Proof}. Write (1.12) as
$$
k(s) := \int_1^\infty E_2^2(x)x^{-s-1}\d x = {2\over s}(
\K(s) + \hf E_2^2(1)),\leqno(3.9)
$$
so that $k(s)$ is regular for $\s>1$. From the Mellin inversion formula
for the (modified) Mellin transform (see [7, Lemma 1]) we have
$$
E_2^2(x) = {1\over2\pi i}\int_{(c)}k(s)x^s\d s\qquad(x >1,\,c>1).
$$
If $\psi(t)$ is a smooth, nonnegative function supported in
$[T/2,\,5T/2]$ such that $\psi(t) = 1$ for $T \le t \le 2T$, then
$$
\int_T^{2T} E_2^4(x)\d x \le
\int_{T/2}^{5T/2}\psi(x) E_2^4(x)\d x
= {1\over2\pi i}\int\limits_{(c)}k(s)
\Bigl(\int_{T/2}^{5T/2}\psi(x) E_2^2(x)
 x^s\d x\Bigr)\d s.\leqno(3.10)
$$
In the last integral over $x$ we perform a large number of integrations
by parts, keeping in mind that $\psi^{(j)}(x) \ll_j T^{-j}\;(j = 0,1,
\ldots\,)$. It transpries that only the values of $|t| \le T^{1+\e}$
in the integral over $s = \s+it$ will make a non-negligible contribution.
Hence (3.10) (with $c =\s > 1$) and Lemma 1 yield
$$\eqalign{I &:=
\int_{T/2}^{5T/2}\psi(x) E_2^4(x)\d x
\ll_\e 1 + \int_{-T^{1+\e}}^{T^{1+\e}}|k(\s+it)|
\left|\int_{T/2}^{5T/2}\psi(x) E_2^2(x)x^s\d x\right|\d t\cr&
\ll_\e 1 + \left(\int_{-T^{1+\e}}^{T^{1+\e}}|k(\s+it)|^2\d t\right)^{1/2}
\left(\int_{T/2}^{5T/2}\psi^2(x) E_2^4(x)x^{2\s+1}\d x\right)^{1/2}\cr&
\ll_\e 1 + \left(\int_{0}^{T^{1+\e}}|k(\s+it)|^2\d t\right)^{1/2}
T^{\s+{1\over2}}I^{1/2}.
\cr}
$$
Simplifying the above expression and using (3.9) we arrive at (3.8).

\medskip
One expects, in conjunction with the conjecture $E_2(x) \ll_\e x^{1/2
+\e}$, the bound
$$
\int_0^T E_2^4(t)\d t \ll_\e T^{3+\e}\leqno(3.11)
$$
to hold as well. In fact, from the author's work [5, Theorem 2] with
$a = 4$, one sees that the lower bound
$$
\int_0^T E_2^4(t)\d t \gg T^{3}
$$
does indeed hold. The upper bound in (3.11) nevertheless seems
unattainable at present. If true, it implies (by e.g., [7, eq. (4.4)]
and H\"older's inequality) the hitherto unproved bounds 
$E_2(T) \ll_\e T^{3/5+\e}$ and ([2, Lemma 4.1]) $\zt \ll_\e |t|^{3/20+\e}$.
From (3.2) of Theorem 2 with $\s = 7/6$ we obtain
$$
\int_0^T E_2^4(t)\d t \ll_\e T^{10/3+\e}.\leqno(3.12)
$$
However, the bound (3.12) was already used in proving Theorem 3 via (3.7).
It is (up to `$\e$') the strongest known bound for the integral in question.
From (3.8) it is seen that the conjectural bound (3.11) holds if
$$
\int_1^T|{\cal K}(\s+it)|^2\d t \ll_\e  T^{2+\e}\qquad(\s > 1)\leqno(3.13)
$$
holds. Conversely, if (3.11) holds, then the bound in (3.7) is to be replaced
by $\ll_\e T^2X^{2-2\s+\e}\;(\s>1)$, and (3.13) follows from this bound
and (3.6) (with $X=T^{6/5}$). Therefore (3.11) is equivalent to the mean
square bound (3.13).

\bigskip
{\hh 4. The function $J(s,\xi)$}
\bigskip
The result on the function $J(s,\xi)\;(0 \le \xi < 1)$ is contained in

\medskip
THEOREM 4. {\it The function $J(s,\xi)$ admits analytic continuation
to the region $\R s > \hf$, where it represents a regular function.
Moreover}
$$
J(\s+it,\xi) \ll_\e t^{-1} + t^{{1-{1\over2}\xi - \s\over1-\xi}+\e}
\qquad(\s > \hf, \;t \ge t_0, \;0 \le \xi < 1).\leqno(4.1)
$$

\bigskip {\bf Proof}.
Let $X = t^{1/(1-\xi)-\delta}$ for a small, fixed $\delta >0$.
We define a sequence of non-negative, smooth 
functions $\rho_j(x)\;(j\in\NN)$  in the following way.
Let $\rho_1(x) \,(\ge 0)$ be a smooth function supported in
$[1, 2X]$ such that $\rho_1(x) = 1$ for $1\le x \le X$, and
$\rho_1(x)$ monotonically decreases from 1 to 0 in $\,[X,\,2X]\,$.
The function $\rho_2(x)$ is  supported in $[X,\,6X]\,$,
where $\rho_2(x) = 1 - \rho_1(x)$ for $X \le x \le 2X$,
$\rho_2(x) = 1$ for $2X \le x \le 4X$  and 
$\rho_2(x)$ monotonically decreases from 1 to 0 in $\,[4X,\,6X]\,$.
In general,  the function $\rho_j(x)$,
supported in $[2^{j-1}X,\,3\cdot2^{j}X]$, satisfies 
$\rho_j(x) = 1 - \rho_{j-1}(x)$ for $2^{j-1}X \le  x \le 3\cdot2^{j-1}X$, 
$\rho_j(x) = 1$ for $[2^{j-1}X, 2^jX]$ and then decreases monotonically
from 1 to 0 in $[2^jX, \,3\cdot2^{j}X]$.
In this way we obtain that 
$$
\rho_j^{(r)}(x) \;\ll_{j,r}\; (2^jX)^{-r}\qquad(j,r \in \NN).\leqno(4.2)
$$
Now we write (cf. (1.16))
$$
J(s,\xi) = \int_1^{2X}\rho_1(x)I(x;\xi)x^{-s}\d x
+ \sum_{j\ge 2}\int_{2^{j-1}X}^{3\cdot2^{j}X}\rho_j(x)I(x;\xi)x^{-s}\d x.
\leqno(4.3)
$$
In the first integral in (4.3) we insert the expression (1.16) 
for $I(x;\xi)$ and integrate repeatedly by parts the factor $x^{-1/2-
i\k_j}$. The integrated terms, after $r$ integrations by parts, will be
$$
\sum_{j=1}^r {A_j\over (s-\hf)^j}
$$
for suitable constants $A_j$. The remaining integral, in view of (4.2),
will be $\ll t^{-B}$ for any given $B>0$, provided that $r = r(B)$ is
sufficiently large. There remain the integrals 
$$
I(K) := \int_{K/2}^{3K}\rho(x)I(x;\xi)x^{-s}\d x
\qquad(\rho(x) = \rho_j(x), K = 2^jX).
$$
Writing the sine in (1.16) as a sum of exponentials, it follows that 
$I(K)$ is a linear combination of expressions the type
$$
J_\pm(K) :=  \sum_{j=1}^\infty \a_j H_j^3(\hf)\k_j^{-1/2}
{\rm e}^{\pm\k_j\log{\k_j\over4{\rm e}}}
\int_{K/2}^{3K}\rho(x)x^{-{1\over2}\mp i\k_j-s}
\exp(-{\txt{1\over4}}(x^{\xi-1}\k_j)^2) \d x,
$$
and we may consider only the case of the `+' sign, since the other 
case is treated analogously. The above series may be truncated at $\k_j
= K^{1-\xi}\log K$ with a negligible error. After an integration by
parts the  integral in $J_\pm(K)$ becomes
$$\eqalign{&
{1\over s - i\k_j+\hf}\int_{K/2}^{3K}x^{{1\over2}+ i\k_j-s}
\exp(-{\txt{1\over4}}(x^{\xi-1}\k_j)^2)\times\cr&
\times\left(\rho'(x) + \hf(1-\xi)\rho(x)x^{2\xi-3}\k_j^2\right)\d x.\cr}
$$
In the range $\k_j \le K^{1-\xi}\log K$ the above expression in
parentheses is
$$
\ll K^{-1} + K^{2\xi-3}K^{2-2\xi}\log^2K \ll K^{-1}\log^2K.
$$
It transpires that, performing sufficiently many integrations by parts,
only the values of $\k_j$ for which $|\k_j - t| \le K^\e$ will make
a non-negligible contribution.
For the estimation of $\a_j\H$ in short intervals we shall need (see
the author's work [6])

\medskip
LEMMA 2. {\it We have}
$$
\sum_{K-G\le\k_j\le K+G} \a_j\H \;\ll_\e\; GK^{1+\e}
\quad(K^{\e}  \;\le\; G \;\le \; K).\leqno(4.4)
$$

\medskip\no
Note that (4.4) implies the bound
$$
H_j(\hf) \;\ll_\e\; \k_j^{1/3+\e},
$$
which breaks the convexity bound $H_j(\hf) \;\ll_\e\; \k_j^{1/2+\e}$, 
but is still far away from the conjectural bound
$$
H_j(\hf + it) \;\ll_\e\; (\k_j + |t|)^\e,
$$
which may be thought of as the analogue of the classical 
Lindel\"of hypothesis ($\zt \ll_\e |t|^\e$) for the Hecke series. 

\medskip
To complete the proof of Theorem 4, note that with the use
of Lemma 2 we obtain
$$\eqalign{
J_\pm(K) & \ll_\e K^{-1/2-\s}K\sum_{|\k_j-t|\le K^\e}\a_j\H \k_j^{-1/2}\cr&
\ll_\e K^{1/2-\s}t^{1/2+\e} \ll_\e 
t^{{1-{1\over2}\xi - \s\over1-\xi}+\e}
\cr}
$$
since $K \gg X (= t^{1/(1-\xi)-\delta})$. This leads to (4.1)
 in view of (4.3) and the preceding
discussion.

\vfill
\eject

\centerline{\hh REFERENCES}

\item {[1]} A. Ivi\'c,  The Riemann zeta-function, {\it John Wiley and
Sons}, New York, 1985.
\item {[2]} A. Ivi\'c,  Mean values of the Riemann zeta-function, LN's
{\bf 82}, {\it Tata Institute of Fundamental Research}, Bombay,  1991
(distr. by Springer Verlag, Berlin etc.).
\item{[3]}  A. Ivi\'c,  On the fourth moment of the Riemann
zeta-function, {\it Publs. Inst. Math. (Belgrade)} {\bf 57(71)}
(1995), 101-110.
\item{[4]} A. Ivi\'c,  The Mellin transform and the Riemann
zeta-function,  {\it Proceedings of the Conference on Elementary and
Analytic Number Theory  (Vienna, July }18-20, 1996),  Universit\"at
Wien \& Universit\"at f\"ur Bodenkultur, Eds. W.G. Nowak and J.
Schoi{\ss}engeier, Vienna 1996, 112-127.
\item{[5]} A. Ivi\'c, On the error term for the fourth moment of the
Riemann zeta-function, {\it J. London Math. Soc.} {\bf60}(2)(1999), 21-32.
\item{[6]} A. Ivi\'c, On sums of Hecke series in short intervals,
{\it J. de Th\'eorie des Nombres Bordeaux} {\bf13}(2001), 1-16.
\item{[7]} A. Ivi\'c, On some conjectures and results for the
Riemann zeta-function, {\it Acta Arith.} {\bf109}(2001), 115-145. 
\item{[8]} A. Ivi\'c, Some mean value results for the Riemann
zeta-function, in `Number Theory. Proc. Turku Symposium 1999'
(M. Jutila  et al. eds.),  de Gruyter, 2001, Berlin, 145-161.
\item{[9]} A. Ivi\'c, On the estimation of ${\cal Z}_2(s)$, in
`Anal. Probab. Number Theory' (A. Dubickas et al. eds.),
TEV, 2002, Vilnius, 83-98.
\item{ [10]} A. Ivi\'c and Y. Motohashi,  The mean square of the error
term for the fourth moment of the zeta-function,  {\it Proc. London
Math. Soc.} (3){\bf 66}(1994), 309-329.
\item {[11]} A. Ivi\'c and Y. Motohashi,  The fourth moment of the
Riemann zeta-function, {\it J. Number Theory} {\bf 51}(1995), 16-45.
\item {[12]} A. Ivi\'c, M. Jutila and Y. Motohashi, The Mellin
transform of powers of  the Riemann zeta-function,  {\it Acta Arith.}
{\bf95}(2000), 305-342.
\item{ [13]} Y. Motohashi,   An explicit formula for the fourth power
mean of the Riemann zeta-function, {\it Acta Math. }{\bf 170}(1993),
181-220.
\item {[14]} Y. Motohashi,  A relation  between the Riemann
zeta-function and the hyperbolic Laplacian, {\it Annali Scuola Norm.
Sup. Pisa, Cl. Sci. IV ser.} {\bf 22}(1995), 299-313.
\item {[15]} Y. Motohashi,  Spectral theory of the Riemann
zeta-function, {\it Cambridge University Press}, Cambridge, 1997.
\medskip
\parindent=0pt

\cc
Aleksandar Ivi\'c \par
Katedra Matematike RGF-a\par
Universitet u Beogradu\par
\DJ u\v sina 7, 11000 Beograd\par
Serbia and Montenegro,
{\sevenbf
aivic@rgf.bg.ac.yu}

\bye